\documentclass[11pt]{article}
\usepackage{latexsym}
\usepackage{amsmath}
\usepackage{amssymb}
\newenvironment{proof}{\noindent {\em Proof. }}{\hfill $\Box$}
\newtheorem {nummer} {} [section]

\newcommand {\A} {{\rm Alt}}

\renewcommand{\dim}  {\mathrm{dim}}

\newcommand{\GL} {{\rm GL}}

\renewcommand {\phi} {\varphi}

\newcommand{\Ref}[1]{\mbox{(}\ref{#1}\mbox{)}}

\newcommand{\SL} {{\rm SL}}

\newcommand{\Sp} [1] {\mbox{{\rm Sp}$(#1)$}}
\newcommand{\FSp} [1] {\mbox{{\rm FSp}$(#1)$}}

\newcommand {\F} {\mathbb{F}}

\newcommand {\PP} {\mathbb{P}}

\def\cc{{\mathchoice
   {\setbox0=\hbox{$\displaystyle\rm C$}\hbox{\raise 0.15
      \ht0\hbox to0pt{\kern0.4\wd0\vrule height0.8\ht0\hss}\box0}}
   {\setbox0=\hbox{$\textstyle\rm C$}\hbox{\raise 0.15
      \ht0\hbox to0pt{\kern0.4\wd0\vrule height0.8\ht0\hss}\box0}}
   {\setbox0=\hbox{$\scriptstyle\rm C$}\hbox{\raise 0.15
      \ht0\hbox to0pt{\kern0.4\wd0\vrule height0.7\ht0\hss}\box0}}
   {\setbox0=\hbox{$\scriptscriptstyle\rm C$}\hbox{\raise 0.15
      \ht0\hbox to0pt{\kern0.4\wd0\vrule height0.7\ht0\hss}\box0}}
}}

\usepackage{fancyhdr}
\begin{document}
\title {On some subgroups of linear groups over $\F_2$\\
 generated by elements of order $3$}
\author{Hans Cuypers\\
email: f.g.m.t.cuypers@tue.nl}
\date{\today}
\maketitle

\bigskip

\begin{abstract}
\noindent
Let $V$ be a vector space over the field $\F_2$.
We investigate subgroups of the linear group
$\GL(V)$ which are generated by a conjugacy class $D$ of
elements of order $3$ such that  all $d\in D$
have $2$-dimensional commutator space $[V,d]$.\\

\noindent
MSC 2010: 20D06, 20E45, 51A45, 51A50 \\
Keywords: Linear groups, elements of order $3$  
\end{abstract}

\newpage

\section{Introduction}

In his revision of Quadratic Pairs \cite{chermak1,chermak2}, 
Chermak \cite{chermak1}
classifies  various subgroups of the symplectic groups $\Sp{2n,2}$
generated by elements $d$ of order $3$ with $[V,d]=\{vd-v\mid v\in V\}$ 
being $2$-dimensional,
where $V$ is the natural module of $\Sp{2n,2}$.
Besides the full symplectic group he encounters orthogonal and unitary groups
over the field with $2$ or $4$ elements, respectively,
as well as alternating groups.
Chermak's proof of this classification
is inductive and relies mainly on methods from geometric algebra. 

By using discrete 
geometric methods we are able to classify
subgroups of ${\rm GL}(V)$, where $V$ is an $\F_2$-vector space
of possibly infinite dimension, generated by 
elements $d$ of order $3$ with $[V,d]$ being $2$-dimensional.

In particular, we prove the following.
(For notation and definitions, the reader is referred to the next section.)

\begin{nummer}{\bf Theorem.}\label{main}
Let $V$ be a vector space of dimension at least $3$ 
over the field $\F_2$. Suppose $G\leq  \GL{(V)}$ is a group
generated by a conjugacy class $D$ of elements of order $3$ such that
\begin{enumerate}
\item $[V,d]$ is $2$-dimensional for all $d\in D$;
\item $[V,G]=V$ and $C_V(G)=\{0\}$.
\end{enumerate} 
Then, up to isomorphism, we have one of the following.
\begin{enumerate}
\item There exists a subspace $\Phi$ of $V^*$ annihilating $V$ such that
$G={\rm T}(V,\Phi)$;
the class $D$ is the unique class of elements of order $3$ with
$2$-dimensional
commutator on $V$.
\item $\dim(V)=3$ and $G= 7:3$; $D$ is one of the two classes of elements of
  order $3$ in $G$.
\item $\dim(V)=4$ and $G={\A}_7$ (inside $\A_8\simeq {\rm GL}(4,2)$); the class $D$ corresponds to the
  class of elements of order $3$ which are products of two disjoint $3$-cycles
  inside $\A_7$.
\item $\dim(V)\geq 6$, and 
$G=\FSp{V,f}$ with respect to some nondegenerate symplectic form $f$ on $V$;
the class $D$ is the unique class of elements of order $3$ with
$2$-dimensional
commutator on $V$.
\item $\dim(V)\geq 6$ and 
$G={\rm F}\Omega(V,Q)$ for some nondegenerate quadratic form $Q$ on $V$ with trivial radical;
the class $D$ is the unique class of elements of order $3$ with
$2$-dimensional
commutator on $V$.
\item $G={\rm FAlt}(\Omega)$ for some set $\Omega$
of size at least $5$;
the class $D$ corresponds to the  class of $3$-cycles, or in case
$|\Omega|=6$,
the class of elements which are a  products of two disjoint $3$-cycles.
The space $V$ is the subspace of the space $\F_2\Omega$ generated by
all vectors of even weight, or, in case $|\Omega|$ is even, the quotient of this subspace
by the all one vector.
\item  $V$ carries a $G$-invariant structure 
of an $\F_4$-space $V_4$ such that $G$ equals ${\rm R}(V_4,\Phi)$, where $\Phi$ is a subspace of $V_4^*$ annihilating
$V_4$. The class $D$ is the class of reflections in $G$.
\item $V$ carries a $G$-invariant structure $(V_4,h)$ 
of an $\F_4$-space $V_4$ equipped
with a  nondegenerate Hermitian  form $h$
such that $G={\rm FU}(V_4,h)$,
the subgroup of ${\rm GU}(V_4,h)$ generated by all reflections.
The class $D$ is the class of reflections in $G$.
\end{enumerate}
\end{nummer}

As indicated above, our proof of this theorem is of  geometric nature.
By using methods similar to those developed in 
Cameron and Hall \cite{ch} and  Cohen, Cuypers and Sterk \cite{ccs}
we are able to show that the subspaces $[V,d]$ with $d\in D$,
are  either all the lines of $\PP(V)$ (leading to the cases (a)-(c) of the
theorem) 
or of a cotriangular space embedded
in $\PP(V)$ (cases (d)-(f)) or these spaces are  (part of) the  one-dimensional subspaces
of an $\F_4$-space induced on $V$ (leading to the cases (g) and (h)).
Although not entirely self contained (we rely on Jonathan Hall's
classification
of cotriangular spaces \cite{hall}), our proofs and  methods are completely elementary.

In the following section we describe the examples occurring in the conclusion
of Theorem \ref{main} somewhat closer. 
The Sections 3 and 4 are devoted to the proof of Theorem \ref{main}.
In particular, in Section 3 
we consider the case where there are $d,e\in D$ with $[V,d]\cap [V,e]$ being
one-dimensional, leading to the groups defined over $\F_2$,
while the final section covers the remaining cases of groups defined over $\F_4$.

\section{The examples and their geometries}
\label{examples}
Suppose $V$ is a vector space over the field $\F_2$.
If $\dim(V)< \infty$, then the
generic example of a group $G$ generated by a class $D$ of elements
$d$ with $[V,d]$ of dimension $2$ is the 
group $\SL(V)$. The elements in $D$  correspond to
conjugates of the element 

$$\left(
\begin{array}{cccccc}
0&1&0&0&\dots&0\\
1&1&0&0&\dots&0\\
0&0&1&0&\dots&0\\
0&0&0&1&     &0\\
\vdots& & & &\ddots & \vdots\\
0&0&0&0&\dots &1\\
 \end{array} 
\right).$$

In small dimensions we encounter two exceptional examples
of groups satisfying the hypothesis of Theorem \ref{main}.

Suppose $V$ has dimension $4$. The group 
$\SL(V)$ is isomorphic to $\A_8$. Under this
isomorphism the elements
in $D$ correspond to those elements in $\A_8$ that are products of two
disjoint $3$-cycles.
The subgroup $\A_7$ of $\A_8$ also acts irreducibly on $V$ and 
is  of course generated by its elements of $D$.

Suppose $\dim(V)=3$. Then the group $\SL(V)$ contains a 
subgroup $7:3$, a split extension of a group of order $7$ by a group of
order $3$, which is irreducible on $V$. This subgroup is
generated by its elements of order $3$.

If $V$ is infinite dimensional, then 
the set of all elements of $\GL(V)$ of order $3$
with two dimensional commutator generates the subgroup ${\rm FSL}(V)$
of $\GL(V)$, 
consisting of all the finitary elements of determinant $1$, see \cite{ch}.
Inside this group we encounter more examples which we now describe.
Let $0\neq v\in V$ and $0\neq \phi\in V^*$ with $\phi(v)=0$.
Then $t_{v,\phi}$ denotes the \emph{transvection}
$$t_{v,\phi}:V\rightarrow V, w\in V\mapsto w+\phi(w)v.$$
The group $G={\rm T}(V,\Phi)$ where $\Phi$ is a subspace of $V^*$
is defined to be the subgroup of ${\rm GL}(V)$ generated by
all transvections $t_{v,\phi}$ with $0\neq v\in V$ and  $0\neq \phi\in \Phi$ with
$\phi(v)=0$. Suppose $\Phi$ annihilates $V$. Then one easily checks that
 $C_V(G)=\{0\}$.
If $v,w\in V$ and $\phi,\psi\in \Phi$ such that
$\phi(v)=\psi(w)=0$ but $\psi(v)=\phi(w)=1$, then 
the product 
$d=t_{v,\phi}t_{w,\psi}$
 is an element of order $3$ with commutator $[V,d]=\langle v,w\rangle$
of dimension $2$.
Let $D$ denote the set of all such elements. It is straightforward to 
check that this set $D$ is a conjugacy class of $G$ generating $G$.
Clearly, $[V,d]$ with $d$ running through $D$ is the set of all
$2$-spaces of $V$. So, $[V,G]=V$.

Next suppose that $(V,f)$ is a 
symplectic space over $\F_2$. The \emph{radical} $\mathrm{Rad}(f)$ is the
subspace $\{v\in V| f(v,w)=0$ for all $w\in V\}$ of $V$. 
For any nonzero vector $v\in V\setminus \mathrm{Rad}(f)$  the transvection
$$t_v:V\rightarrow V, w \in V\mapsto w+f(v,w)v$$
is a nontrivial element of ${\rm Sp}(V,f)$.
It is well known that the set of such transvections forms a
conjugacy class of $3$-transpositions in ${\rm Sp}(V,f)$.
If we set $D$ to be the set of all products $t_vt_w$
where $v,w\in V$ with $f(v,w)=1$, then one readily checks that
$D$ is a conjugacy class
of elements of order $3$ in $G=\langle D\rangle$ with $[V,d]$ of dimension $2$.
The radical $\mathrm{Rad}(f)$ is contained in $C_V(G)$.

By ${\cal S}p(V,f)$
we denote the partial linear 
space  $(P,L)$ where $P$ consists of all  the nonzero vectors
of  $V$.
A line in $L$ is  the set of three nonzero vectors in a $2$-dimensional 
subspace $W$ of $V$ on which $f$ does not vanish.

Notice that for each $d\in D$, the commutator $[V,d]$ determines a unique line
of $\mathcal{S}p(V,f)$.
The class $D$ generates $\FSp{V,f}$
provided $\dim(V)\geq 6$. Here $\FSp{V,f}$ 
denotes the finitary symplectic group, i.e., the subgroup of finitary elements
in the symplectic group $\Sp{V,f}$.
If $\dim(V)=4$ or $2$, then $D$ generates the derived groups
$\A_6\leq  {\rm Sp}(4,2)$ and $\A_3\leq  {\rm Sp}(2,2)$, respectively.

Next consider 
a quadratic form $Q$ on $V$ whose associated bilinear form is the
symplectic form $f$.
The \emph{radical} ${\rm Rad}(Q)$ of $Q$ is defined to be the radical of $f$.
The  transvection $t_v$, where 
$v\in V\setminus {\rm Rad}(Q)$ with $Q(v)=1$, 
is in the orthogonal group ${\rm O}(V,Q)$.
Suppose $\dim(V)\geq 6$. Then the subset $D_Q$ of $D$ of all
elements obtained as products $t_v t_w$ with 
$Q(v)=Q(w)=1$ and $f(v,w)=1$ is a conjugacy class
of $\Omega(V,Q)$, the derived subgroup of ${\rm O}(V,Q)$.
The set $D_Q$  generates the  finitary group ${\rm F\Omega}(V,Q)$.
The radical of $Q$ is centralized by ${\rm F\Omega}(V,Q)$.

The corresponding geometry $\mathcal{N}(V,Q)$ has as points
the vectors $v\in V$ with $Q(v)=1$.
A typical line is the set of three  nonzero vectors in an elliptic 
$2$-space, i.e., a $2$-space in which $Q(v)=1$ for any nonzero vector $v$
contained in it.
Clearly, ${\cal N}(V,Q)$ is a subspace of ${\cal S}p(V,f)$.

There is yet another class of subgroups of ${\rm Sp}(V,f)$
generated by a subset of $D$. To describe this class we will 
start with a particular description of the  symplectic space $(V,f)$.
Indeed, the symplectic space $(V,f)$ might be obtained as follows.
Suppose $\Omega$ is 
a (possibly infinite) set. 
Let $\F_2\Omega$ be the $\F_2$-vector space with basis $\Omega$.
By $E\F_2\Omega$ we denote the subspace
of $\F_2\Omega$  generated by the vectors $\omega_1+\omega_2$, where
$\omega_1,\omega_2\in \Omega$.
Notice that the standard dot product on $\F_2\Omega$ induces a symplectic
form on $E\F_2\Omega$. We can identify $(V,f)$ with this symplectic space.
The transpositions in ${\rm Sym}(\Omega)$ induce transvections
on $V$. So, the $3$-cycles in ${\rm Sym}(\Omega)$, 
(i.e., the products of two noncommuting transpositions)
induce a subset $D_\Omega$ of $D$ generating a subgroup
of ${\rm GL}(V)$ isomorphic to the (finitary) alternating group ${\rm
  FAlt}(\Omega)$, the subgroup of ${\rm Sym}(\Omega)$ of all
permutations with finite support that are even.
The corresponding geometry $\mathcal{T}(\Omega)$
has as points the  vectors 
 $\omega_1+\omega_2$, where
$\omega_1\neq \omega_2\in \Omega$, a line being  the triples of
points of the form  $\omega+\omega'$, where $\omega\neq\omega'$ are taken 
from some subset of size $3$ of $\Omega$.

The above geometries are all examples of {\em cotriangular spaces}.
These are partial linear spaces with lines of size three and having 
the property that any point $p$ not on a line 
$l$ is collinear with no point or with exactly two
points on $l$. A cotriangular space is called \emph{irreducible}
if it is connected and for any pair of points
$p,q$ we have that $p^\perp=q^\perp$ implies $p=q$.
Here $p^\perp$ denotes the set consisting of $p$ and
all points \emph{not} collinear with $p$.
The spaces described above
are characterized by the following result (rephrased to fit our purposes) 
of Jonathan Hall.

\begin{nummer}{\rm \bf (J.I. Hall \cite{hall})}\label{cotrianglethm}
Let $V$ be a vector space over the field $\F_2$.
 Let $\Pi=(P,L)$ be an irreducible cotriangular space, where $P$
is a subset of $V\setminus \{0\}$, and each line in $L$ is a triple of points
inside a $2$-dimensional subspace of $V$.
If $P$ generates $V$ and $\bigcap_{p\in P}\langle p^\perp\rangle=\{0\}$
we have one of the following:
\begin{enumerate}
\item $\Pi=\mathcal{S}p(V,f)$ for some nondegenerate symplectic form $f$ 
on $V$.
\item $\Pi=\mathcal{N}(V,Q)$ for some nondegenerate quadratic form $Q$ on $V$
with trivial radical.
\item There is a set $\Omega$ such that
$\Pi=\mathcal{T}(\Omega)$ and  $V=E\F_2(\Omega)$ or, 
in case $|\Omega|$ is even, $V=E\F_2(\Omega)/\langle\Sigma_{\omega\in \Omega}\ \omega\rangle$, the 
quotient of $E\F_2(\Omega)$ by the all one vector.
\end{enumerate}

\end{nummer}

Finally we shall discuss the examples 
coming from groups defined over $\F_4$.
Let $V_4$ be a vector space over $\F_4$.
For every $v\in V_4$ and $\phi \in V_4^*$, with $\phi(v)\neq 0,1$ 
we define 
the map $$r_{v,\phi}: V\rightarrow V_4, w\in V_4\mapsto w-\phi(w) v.$$
The map $r_{v,\phi}$ is a {\it reflection} with {\it center\/} 
$\langle v\rangle$ and {\it axis\/} $\ker \phi$.
A reflection has order $3$.
If $\Phi$ is a subspace of $V_4^*$, then denote by
${\rm R}(V_4,\Phi)$ the subgroup of ${\rm GL}(V_4)$ generated by all
reflections $r_{v,\phi}$ with $v\in V_4$, $\phi\in \Phi$ and 
$\phi(v)\neq 0,1$.
If $\dim(V_4)$ is finite-dimensional, then
${\rm R}(V_4,V_4^*)={\rm GL}(V_4)$.
Let $V$ denote the space $V_4$ considered as an $\F_2$-space.
The reflections
 provide  examples of elements
of order $3$ having a $2$-dimensional commutator on $V$.

If  $h$ is a nondegenerate Hermitian 
form on $V_4$, then
for each vector $v\in V$ with $h(v,v)=1$ and $\alpha\in \mathbb{F}_4$, $\alpha\neq 0,1$,
the map
$$r_v: w\in V_4 \mapsto w+ \alpha h(w,v)v$$
is a reflection in  the finitary unitary group
$${\rm FU}(V_4,h) = 
\{g\in {\rm FGL}(V)\mid \forall_{x,y\in V} \quad  h(xg,yg)) = h(x,y)\}.$$
In fact, all these reflections
generate
the finitary group ${\rm FU}(V_4,h)$. 

Notice that in these examples over $\F_4$, the commutators
$[V,r_1]$ and $[V,r_2]$, where $r_1$ and $r_2$ are reflections
on $V_4$, either are equal or meet trivially.

\section{Geometries with points and groups over $\F_2$}
\label{points}

Let $V$ be  a vector space over  $\F_2$, 
and suppose $G\leq  {\rm GL}(V)$ is
generated by a normal set $D$ of  elements $d\in G$ of order $3$ 
such that
$[V,d]$ is $2$-dimensional.
(Here normal means closed under conjugation.)

\begin{nummer}\label{directsum}
Suppose $d\in D$. Then $V=[V,d]\oplus C_V(d)$.
\end{nummer}

\begin{proof}
Suppose $v\in V$ with $[v,d]\in [V,d]\cap C_V(d)$.
Then $0=[vd+v,d]=vd^2+vd+vd+v=vd^2+v$. But then
$[v,d]=(vd^2+v)d=0$. We have found that
$[V,d]\cap C_V(d)=0$.
As each $v\in V$ equals 
$v=v+[vd,d]+[vd,d]=v+vd+vd^2+[vd,d]\in C_V(d)+[V,d]$
we have proved that 
$V=[V,d]+C_V(d)$.
\end{proof}

\begin{nummer}\label{invariant}
A subspace $W$ of $V$ is invariant under $d\in D$
if and only if $W\leq  C_V(d)$ or $[V,d]\leq  W$.
\end{nummer}

\begin{proof}
If $W$ is centralized by $d$, it clearly is invariant.
If $[V,d]\leq W$, then for $w\in W$ we have $wd=[w,d]+w\in W$
and we find $W$ to invariant under $d$.

Now suppose $w\in W$ is invariant, but not centralized by $d$. 
Then $0\neq [w,d]\in W\cap [V,d]$.
But then  $[V,d]=\langle [w,d], [w,d]d\rangle$ is contained in $W$.
\end{proof}

\bigskip

A {\em $D$-line}, or just {\em line}, is a subspace of 
$V$ of the form $[V,d]$ with $d\in D$.
A {\em $D$-point}, or just {\em point}, 
is a $1$-space of $V$ which is the intersection of
two distinct $D$-lines.
Both points and lines are also considered to be points and lines of the 
projective space $\mathbb{P}(V)$. 

Let $\mathcal{P}$ be the set of $D$-points and $\mathcal{L}$ the set of
$D$-lines. The geometry $\Pi(D)$ is the pair $(\mathcal{P},\mathcal{L})$,
where incidence is symmetrized containment.
A line is often identified with the set of points it contains.

If $W$ is a subspace of $V$, then by $\Pi(D)_W$ we denote the pair
$(\mathcal{P}_W, \mathcal{L}_W)$ where 
$\mathcal{L}_W$ is the set of  $D$-lines contained in $W$,
and $\mathcal{P}_W$ the set of intersection points
of two distinct lines in $\mathcal{L}_W$ meeting nontrivially.

Let $U$ be a subspace of $V$. Then by $D_U$ we denote the set of all $d\in D$
with $[V,d]\leq  U$. The subspace  $A_U$ of $V$ is equal to $\bigcap_{d\in D_U}\
C_V(d)$.

\begin{nummer}
If $l$ is a $D$-line, then it contains zero or three  $D$-points.
\end{nummer}

\begin{proof}
Since $[V,d]\cap C_V(d)=\{0\}$ by \Ref{directsum}, we find that $d$ is transitive on the three
nonzero vectors in $[V,d]$. This clearly implies the result.
\end{proof}

\begin{nummer}\label{planes}
Suppose $l$ and $m$ are distinct $D$-lines intersecting at a point.
Let $W$ be the subspace $l+m$ of $V$.
Then $\Pi(D)_W$ is either a dual affine plane or a projective plane
in $\PP(W)$.

The group $\langle D_W\rangle$ is transitive on the lines 
in $\mathcal{L}_W$.
\end{nummer}

\begin{proof}
Let $d\in D$ with $[V,d]=l$.
Then $\langle d\rangle$ fixes a unique point in $\mathbb{P}(W)$,
call this point $q$, 
is transitive on the three points of $l$, and transitive on the
three remaining points of $\mathbb{P}(W)$.
The group $\langle d\rangle$ fixes $l$, is transitive on the three
lines of  $\mathbb{P}(W)$ on $q$ and on the three remaining lines.

If $q$ is a point on a $D$-line, then clearly
$\langle D_W\rangle$ is transitive on the points and lines
of $\mathbb{P}(W)$ and  $\Pi(D)_W$ equals $\mathbb{P}(W)$.

If $q$ not on any $D$-line, then   $\Pi(D)_W$ equals
the dual affine plane of all points of $\Pi(D)_W$
different from $q$ and all lines not on $q$.
Also in this case $\langle D_W\rangle$ is transitive on the points and lines
of  $\Pi(D)_W$.
\end{proof}

\begin{nummer}\label{transon}
If $D$ is a conjugacy class in $G$, then $G$ is transitive on $\mathcal{L}$
and $\mathcal{P}$.
\end{nummer}

\begin{proof}
Transitivity of $G$ on $D$ implies transitivity on lines.
As each $d\in D$ is transitive on the three $1$-spaces of   the line
$[V,d]$, transitivity on points follows immediately.
\end{proof}

\begin{nummer}\label{transit}
If $l,m\in \mathcal{L}$ are in the same connected component $\Pi_0$
of $\Pi(D)$
and $l=l_0,\dots,l_k=m$ is a path from $l$ to $m$ inside
$\Pi_0$, with $l_i$ and $l_{i+1}$ intersecting at a point for $0\leq i< k$,
then there is a $g\in \langle D_{l_0},\dots, D_{l_k}\rangle$
with $lg=m$.  
\end{nummer}

\begin{proof}
By \Ref{planes}, there is for $i=0,\dots, k-1$ a $g_i\in \langle D_{l_0},\dots, D_{l_k}\rangle$
with $l_ig_i=l_{i+1}$. But then $g=g_0\cdots g_{k-1}$
maps $l$ to $m$.  
\end{proof}

\begin{nummer}\label{sym}
Suppose $l\neq m\in \mathcal{L}$ are in the same $G$-orbit on $\mathcal{L}$.
Then $l\subseteq A_m$ if and only if $m\subseteq A_l$.
\end{nummer}

\begin{proof}
Suppose $l\subseteq A_m$, then $l$ and therefore also 
$A_l$ is invariant under each
element $d\in D_m$.
So, by \Ref{invariant} we either have $m\subseteq A_l$ or $A_l\subseteq A_m$.
In the latter case the inclusion is
proper since 
$l\in A_m$ but not in $A_l$.

Now suppose $m\not \subseteq A_l$. 
Then $A_l\subset A_m$. Let $g$ be an element in $G$ with $mg=l$, which is the product of a finite number of elements from $D$. Such an element exists by \Ref{transit}. Then 
$A_l=A_mg\subset A_m$.
However, as $g\in G$ is the product of a finite number of elements from $D$,
the subspace $C_V(g)$ has finite codimension in $V$.
As $g$ induces a bijective map 
$A_m/(A_m\cap C_V(g))\rightarrow A_mg/(A_mg\cap C_V(g))$
between finite dimensional spaces, these spaces are equal.
This implies $A_l=A_mg=A_m$ contradicting the above.
\end{proof}

\begin{nummer}\label{connected2}
Suppose $\mathcal{L}$ is a single  $G$-orbit.
Let $d,e\in D$. 
If $\dim(C_V(d)\cap [V,e])=1$, then 
$[V,d]$ and  $[V,e]$ are in the 
same connected component of $\Pi(D)_{[V,d]+[V,e]}$.
\end{nummer}

\begin{proof}
Let $d,e\in D$ with $\dim(C_V(d)\cap [V,e])=1$.
Notice that  $[V,e]\not \subseteq A_{[V,d]}$.
So, by \Ref{sym} we can assume, eventually
after replacing $e$ by some appropriate element in $D_{[V,e]}$, 
that $[V,d]\not\subseteq C_V(e)$.
 
The lines $[V,e]d$ and $[V,e]$ intersect in the point
$C_V(d)\cap [V,e]$ and are contained in the $4$-dimensional
subspace $[V,d]+[V,e]$.
The subspace $W=[V,e]+[V,e]d$ is $3$-dimensional
and $\Pi(D)_W$ is a dual affine or projective plane inside $\mathbb{P}(W)$.

The line $[V,d]$ intersects $\PP(W)$ in a point $p$ of $\PP(W)$.
If this point $p$ is in $\mathcal{P}_W$ we are done.
Thus assume that this point is not in $\mathcal{P}$.
In this case $\Pi(D)_W$ is a dual affine plane 
and $p$ is the unique 1-space of $\PP(W)$ not in this dual affine
plane. In particular, $p$  equals $C_W(e)$.

By the same arguments we can assume that $[V,d]$ meets
$W'=[V,e]+[V,e]d^2$ in the point $q=C_{W'}(e)$ which is the
unique 1-space of $\PP(W')$ not in $\Pi(D)_{W'}$.
As $p\neq q$, we find that  $[V,d]=\langle p,q\rangle\leq  C_V(e)$, which contradicts our assumption.
\end{proof}

\begin{nummer}\label{connected}
Suppose $C_V(G)=\{0\}$ and
$\mathcal{P}\neq \emptyset$.
Then $\Pi$ is connected if and only if $G$ is transitive on $\mathcal{L}$.

Moreover, if $\Pi$ is connected, then its diameter is at most $2$.
\end{nummer}

\begin{proof}
First assume that $\Pi$ is connected. Then by \Ref{transit}
$G$ is transitive  on $\mathcal{L}$. 

Now suppose that $\Pi$ is not connected, but $G$ is transitive on 
$\mathcal{L}$. 
Since $G$ is generated by  $D$,
there are $f,h\in D$ such that $[V,h]$ and $[V,h]f=[V,h^f]$ are in different
components of $\Pi$. Then also   $[V,h]f$ and $[V,h]f^2$  are in different
components.
Since $h^f$ does not centralize both $[V,h]$ and $[V,h]f^2$,
we can assume that there are noncommuting
$d$ and $e$ in $D$ with $[V,d]$ and $[V,e]$ in different
connected components of $\Pi$.  Fix such an element $e$.
By \Ref{connected2} we can assume
that $C_V(e)\cap [V,d]=\{0\}$.

Let $f\in D$ such that $[V,f]$ meets $[V,d]$ in a point and
let $W$ be the subspace $[V,d]+[V,f]$. 
If $\Pi(D)_{W}$  contains
a line $[V,g]$ with $g\in D$  such that $C_V(e)\cap[V,g]$
is $1$-dimensional,
then \Ref{connected2} gives a contradiction with $[V,e]$ being in a 
connected component of $\Pi(D)$ different form the one containing $[V,d]$.
Hence   $\Pi(D)_W$ is a dual affine plane,
moreover, $C_V(e)$ meets $W$ in the unique 1-space $p$ 
not in that dual affine plane.

As $C_V(G)=\{0\}$, there is an element $h\in D$ not centralizing
$p$. But then $C_V(h)$ meets at least one of the lines of the dual
affine plane $\Pi(D)_W$ in a point. Without loss
of generality we may assume this line to be $[V,d]$.
Let $U=[V,d]+[V,h]$.
If $[V,d]$ and $[V,h]$ meet nontrivially,
then $\Pi(D)_{U}$ is a projective plane
containing a line $l\in \mathcal{L}$ which meets $C_V(e)$ in a point.
As above, this leads by \Ref{connected2}   to a contradiction.
Thus $\dim(U)=4$.
But then $U_1=[V,d]+[V,d]h$ and $U_2=[V,d]+[V,d]h^2$
are two distinct $3$-dimensional spaces
on $[V,d]$.
As above, for both $i=1$ or $2$, 
we can assume that $C_V(e)$ meets 
$U_i$  in a 1-space, which is the unique 1-space
of $U_i$  which is not
in $\mathcal{P}_{U_i}$.
But that implies that $C_U(e)=C_U(d)$.

By \Ref{transit} the above reasoning also
applies to $[V,h]$ and $h$, so that
$C_U(e)=C_U(h)$.
But $C_U(h)\neq C_U(d)$, which is a final contradiction.
Hence we have shown that $\Pi$ is connected.

It remains to prove that the diameter of $\Pi$ is at most $2$, provided
$\Pi$ is connected.
So assume $\Pi$ to be conected and let $p,q,r,s$ be 
a path of length $3$ in the collinearity graph
of $\Pi$. Then both $p$ and $s$ have at least $2$ neighbors on the line through $q$ and $r$, see \Ref{planes}. But that implies that they have a common neighbor. So,
indeed, the diameter of $\Pi$ is at most $2$.
\end{proof}

\bigskip
Let $p,q\in \mathcal{P}$ be points. We write $p\sim q$ if
$p$ and $q$ are distinct collinear points of $\Pi$.
By $p\perp q$ we mean that $p$ and $q$ are equal or noncollinear. 
By $p^\sim$ we denote the set of
all points collinear to $p$ (excluding $p$). The complement of $p^\sim$ in $\mathcal{P}$
is the set $p^\perp$.

If for $p$ and $q$ we have $p^\perp=q^\perp$, then
we write $p\equiv q$. The relation $\equiv$ is obviously an equivalence
relation.

\begin{nummer}\label{reduced}
Suppose $\Pi$ is connected.
If $p\neq q\in \mathcal{P}$ with $p\equiv q$, then 
$C_V(G)\cap p+q\neq \{0\}$.
\end{nummer}

\begin{proof}
Suppose $p\neq q\in \mathcal{P}$ with $p\equiv q$.
Notice that $p\perp q$.
Let $r$ be the third point on the projective line through
$p$ and $q$.
If $l\in \mathcal{L}$ is a line on $p$, then $\Pi_W$ is
a dual affine plane, where $W$ is the subspace spanned by $l$ and $q$.
So, each $d\in D_l$ centralizes $r$.

If $p_1$ and $q_1$ are two noncollinear points
in  $\Pi_W$, then the projective line on $p_1$ and $q_1$ contains $r$.
Moreover, if $s\in p_1^\perp$ but not in $q_1^\perp$, then 
either $s$ is collinear to $p$ but not to $q$, or vice verse.
As this contradicts $p\equiv q$, we find that $p_1\equiv q_1$.

But that implies, by connectivity of $\Pi$, that $r$ is
in $C_V(d)$ for each $d\in D$. In particular, $C_V(G)\cap p+q\neq \{0\}$.
\end{proof}

\begin{nummer}\label{radical}
Suppose $G$ is transitive on $\mathcal{L}$, $C_V(G)=\{0\}$, $[V,G]=V$ 
and $\mathcal{P}\neq \emptyset$.

If $(\mathcal{P},\mathcal{L})$ contains no projective planes, then the subspace
$\bigcap_{p\in \mathcal{P}}\ \langle p^\perp\rangle$ of $V$ is equal to
$\{0\}$.
\end{nummer}

\begin{proof}
Suppose $(\mathcal{P},\mathcal{L})$ does not contain projective planes.
Let $p$ be a point in $\mathcal{P}$ and $l$ a line on $p$.
Then there is an element $d\in D$ with $l=[V,d]$.

Let $q$ be a point in $p^\perp$. If $q$ is collinear with 
a point on $l$, then inside the $3$-space $\langle q,l\rangle$ we find a unique $1$-space $r$, which is not in the dual affine plane generated by $q$ and $l$.
This point $r$ is centralized by $d$ and on the line through $p$ and $q$.
But then  $q\in \langle p,r\rangle \subseteq \langle p,C_V(d)\rangle$.

If $q$ is not collinear to a point on $l$, then, by \Ref{connected}, we can find a point $s$
collinear with both $q$ and the point  $p'=pd$ on $l$.
We can assume that $s$ is in $p^\perp$.
(Indeed, if $s$ is collinear to $p$, then we may replace it by the third 
point on  the line through $p'$
and $s$, which is not collinear with $p$.)
Let $t$ be the third point on the line through $s$ and $q$.
Notice that also $t$ is in $p^\perp$ but collinear to $p'$.
By the previous paragraph we find
both $s,t\in \langle p,C_V(d)\rangle$ but then also $q\in \langle p,C_V(d)\rangle$.
This shows that $\langle p^\perp\rangle$ is contained in
$\langle p,C_V(d)\rangle$.
Moreover, as $C_V(d)$ has codimension $2$ in $V$, we find
 $\langle p^\perp\rangle$ to be a proper subspace of $V$ not containing 
$l=[V,d]$.

This implies that the  space 
$\bigcap_{p\in \mathcal{P}}\ \langle p^\perp\rangle$, which is  invariant 
under each
$e\in D$,  does not contain any line from $\mathcal{L}$. 
But then \Ref{invariant} implies
that $\bigcap_{p\in \mathcal{P}}\ \langle p^\perp\rangle$ is centralized 
by each $e\in D$ and hence is contained in $C_V(G)=\{0\}$, 
proving the result.
\end{proof}

\begin{nummer}{\bf Theorem.}\label{geomthm}
Suppose $G$ is transitive on $\mathcal{L}$, $C_V(G)=\{0\}$, $[V,G]=V$ 
and $\mathcal{P}\neq \emptyset$.
Then, up to isomorphism, $\Pi$ is one of the following spaces:
\begin{enumerate}
\item $\mathbb{P}(V)$.
\item $\mathcal{S}p(V,f)$ for some nondegenerate symplectic form $f$ 
on $V$.
\item $\mathcal{N}(V,Q)$ for some nondegenerate quadratic form $Q$ on $V$
with trivial radical.
\item There is a set $\Omega$ such that $V=E\F_2(\Omega)$
or, in case $|\Omega|$ is even, $V=E\F_2(\Omega)/\langle\Sigma_{\omega\in \Omega}\ \omega\rangle$, the quotient by the all one vector,
and $\Pi=\mathcal{T}(\Omega)$.
\end{enumerate}
\end{nummer}

\begin{proof}
By \Ref{planes} and \Ref{connected}, the space $\Pi$ is a connected partial linear space
of order $3$ in which any two intersecting lines generate a subspace isomorphic to a
projective or a dual affine plane.

If all planes are projective, then clearly $\Pi=\mathbb{P}(V)$ and we are in case (a). 
If all planes are dual affine, then by \Ref{radical}
we can apply Theorem \ref{cotrianglethm} and we are in one of the cases
(b), (c) or (d).  

Now to prove the theorem it suffices to show that $\Pi$ cannot contain both a projective and dual affine plane. So, suppose it does.
Let $\pi$ be a projective plane and $p$ a point outside $\pi$.
We claim that $p$ is collinear with all or all but one of the points of $\pi$.
First assume that  $p$ is collinear with some point $q\in \pi$ and let $r$ be the third point on the line
through $p$ and $q$.
If $p^\perp\cap \pi$ is a line, then so is $r^\perp \cap \pi$. This would imply that there is 
a point in $\pi\cap p^\perp\cap r^\perp\subset q^\perp$, which is clearly impossible. So $p^\perp$ meets $\pi$ in at most one point.

Now assume that $\pi \subseteq p^\perp$. As the diameter of $\Pi$ is at most $2$,
there is a point $q$ collinear to $p$ and also to some point in $\pi$.
Since by the above $q^\perp\cap \pi$ contains at most one point,
there is a projective plane $\pi'$ on $q$ meeting $\pi$ in a line.
But then $p^\perp$ meets $\pi'$ in a line, which contradicts the above, and  we have proved our claim.

Now suppose $p$ and $q$ are noncollinear points.
By \Ref{transon} each line on $p$ is in a projective plane
and hence contains at most one point in $q^\perp$.  
So, all points collinear to $p$ are also collinear to $q$. Similarly
 all points collinear to $q$ are also collinear $p$. This implies
that $p\equiv q$, and,  by \Ref{reduced}, contradicts that $C_V(G)=\{0\}$.
\end{proof}

\begin{nummer}{\bf Theorem.}\label{linear}
Suppose  $\Pi=\mathbb{P}(V)$.
Then, up to isomorphism, $G={\rm T}(V,\Phi)$ for some subspace $\Phi$ of $V^*$
annihilating $V$, or $\dim(V)=4$ and $G=\A_7$,
or $\dim(V)=3$ and $G=7:3$. 
\end{nummer}

\begin{proof}
First assume that the group $G$ contains a transvection $\tau$.
As $G$ is transitive on the points in $\mathbb{P}(V)$, see \Ref{transon},
each point in $\mathbb{P}(V)$ serves as center of some transvection 
in $G$. Suppose $H$ is a hyperplane of $V$ serving as the axis of some
transvection $\tau\in G$. Let $p$ be the center of this transvection.
If $q$ is now a second point in $H$, then let $e$ be an element of $D$
with $p,q\subseteq [V,e]$. Then $q=pe$ or $q=pe^{-1}$ and 
the transvection with center $q$ and axis $H$ is a conjugate of $\tau$ in $G$.

Now let $K$ be a second hyperplane of $\PP(V)$
serving as transvection axis for some
transvection in $G$. Then by the above we can find transvections
$\tau $ and $\sigma$ in $G$ with the same center and with axis $H$ and
$K$ respectively. 
But then $\sigma\tau$ is a transvection with axis the unique hyperplane
$L$ distinct from $H$ and $K$ containing $H\cap K$.
So the elements of $V^*$ serving as  transvection axis for some
transvection in $G$ form the set of nonzero vectors of a subspace $\Phi$ of $V^*$.
This implies that the transvections in $G$ generate the subgroup
${\rm T}(V,\Phi)$ of $G$.

If $\dim(V)=3$ or $4$, and $G$ contains a transvection,
then by the above $G=\SL(V)$. 
So, assume that $G$ does not contain any transvection.

If $\dim(V)=3$, then any involution in $\GL(V)$ is a transvection.
So, $|G|$ is odd and $|G|\mid 21$. On the other hand, $G$ is transitive
on the $7$ lines, while an element $d\in D$ fixes a line.
Hence $G$ has order $21$ and is isomorphic to $7:3$.

If $\dim(V)=4$, then $G$ has order divisible by 
$3\cdot 15\cdot 7$ as the stabilizer of a point-line flag has order at least $3$.
Indeed, an element  $d\in D$, fixes a point-line-flag.
An easy computation within $\A_8\simeq \GL(4,2)$ reveals that 
$G\simeq \A_7$.

Now assume that $\dim(V)\geq 5$.
Fix an element $d\in D$ and the consider the line $[V,d]$.
This line is contained in $5$-dimensional subspace
$U$ of $V$.
The intersection  of  $U$ with $C_V(d)$ is $3$-dimensional.
Pick two lines $l$ and $m$ in $\mathcal{L}$ spanning $C_V(d)\cap\Delta$.
The above shows that 
inside the subgroups generated by $D_{[V,d]+l}$ and  $D_{[V,d]+m}$, 
respectively,
we can find  elements $e\in D_l$ and $f\in D_m$
not  centralizing  $[V,d]$.
But then it is straightforward to check that among the conjugates
of $d$ under $\langle e,f\rangle$ we find two elements,
$d_1$ and $d_2$ say, with $[V,d_1]$ and $[V,d_2]$ meeting at a point.
Moreover, as both $e$ and $f$ leave $C_V(d)$ invariant, we have
$C_V(d_1)=C_V(d_2)=C_V(d)$.
But then either $d_1d_2$ or $d_1d_2^{-1}$ induces a transvection on
$[V,d_1]+[V,d_2]$ with center
$[V,d_1]\cap [V,d_2]$. 
But as $C_V(d)$ is centralized by  $d_1d_2$ or $d_1d_2^{-1}$, 
we have found a transvection $\tau$ on $V$ in $G$.
Now notice that $d_1\in \langle \tau, \tau^{d_1}\rangle$.
So $G$ is generated by its transvections. By the above we can conclude that
$G$ equals ${\rm T}(V,\Phi)$, where $\Phi$ is some subspace of $V^*$.
Since $\bigcap_{\phi\in \Phi}\ker{\phi}$ is centralized by $G$, 
we can conclude that $\bigcap_{\phi\in \Phi}\ker{\phi}=\{0\}$
and $\Phi$ annihilates $V$.
\end{proof}

\begin{nummer}{\bf Theorem.}\label{cotriangle}
Suppose $\Pi$ is a nondegenerate cotriangular space as in case (b), (c) or (d) of \Ref{geomthm}.
Then, up to isomorphism, we have one of the following.

\begin{enumerate}
\item $G= {\rm FSp}(V,f)$ for some nondegenerate symplectic form $f$.
\item $G={\rm F}\Omega(V,Q)$ for some nondegenerate quadratic form $Q$ with trivial radical.
\item $G={\rm FAlt}(\Omega)$ for some set $\Omega$, where  $V=E\F_2\Omega$
or, in case $|\Omega|$ is even, $V=E\F_2(\Omega)/\langle\Sigma_{\omega\in \Omega}\ \omega\rangle$, the quotient by the all one vector.
\end{enumerate}
In all cases $D$ is uniquely determined.
\end{nummer}

\begin{proof}
For each line $l$ of $\Pi$, there is (up to taking inverses) 
at most one element $d\in \GL(V)$ with $[V,d]=l$ and
centralizing the codimension $2$ 
subspace $\bigcap_{p\in l} \langle p^\perp\rangle$ of $V$.
So, this element is in $D$.
But now it is straightforward to check that the theorem holds.
\end{proof}

\bigskip

The above results  classify all the groups satisfying the hypothesis
of Theorem \ref{main} for which the set  
$\mathcal{P}$ is nonempty.

\section{Pointless geometries and groups over $\F_4$}

Suppose $V$ is an $\F_2$-vector space and $G\leq  \GL(V)$ a subgroup
generated by a normal set  $D$ of elements $d\in G$ of order $3$
with $[V,d]$ of dimension $2$.
We keep the notation of the previous section.

In this final section we consider the case where $\mathcal{P}$ is the empty
set.
Although the set $\mathcal{P}$ is empty, 
we will still be able to construct a useful geometry.
However, now the elements of $\mathcal{L}$ will play the role of `points'
and certain $4$-dimensional subspaces of $V$ will play the role
of `lines'. We make this precise in the sequel of this section.

Assume throughout this section that the set $\mathcal{P}$ is empty.
 
\begin{nummer}\label{intersection}
If $d,e\in D$, then $[V,d]\cap C_V(e)=\{0\}$ or $[V,d]\leq  C_V(e)$.
\end{nummer}

\begin{proof}
If  $[V,d]\cap C_v(e)\neq \{0\}$, then $[V,d]e=[V,d^e]$ meets $[V,d]$ nontrivially.
By the assumption that $\mathcal{P}$ is empty, we find 
$[V,d]e=[V,d]$.
By \Ref{directsum} and \Ref{invariant} we find $[V,d]\leq  C_V(e)$. 
\end{proof}

\bigskip

A {\em spread} of a $4$-dimensional subspace $W$ of $V$ is a set of $5$
subspaces of $W$ of dimension $2$, pairwise intersecting in $\{0\}$.

\begin{nummer}\label{defspread}
Suppose $d,e\in D$ with $[V,d]\neq [V,e]$.
If $[V,e]\not \subseteq C_V(d)$, then 
$W:=[V,\langle d,e\rangle]$ is a $4$-dimensional subspace of $V$ containing
$4$ or $5$ lines from $\mathcal{L}$.

If $W$ contains $5$ lines of $\mathcal{L}$, then these $5$
lines form a spread in $W$.

If $W$ contains exactly $4$ lines from $\mathcal{L}$, 
then these $4$ lines together
with the $2$-dimensional space $C_V(d)\cap W$ form a spread of $W$. 
\end{nummer}

\begin{proof}
Let $d,e$ be  elements in $D$ with $[V,d]\neq [V,e]$. 
Then, as by assumption
$[V,d]\cap [V,e]=\{0\}$,  the space 
$W:=[V,d]+[V,e]$ is $4$-dimensional.

If $d$ does not centralize $[V,e]$, then
$W$ contains the four lines $[V,d]$, $[V,e]$, $[V,e]d$ and $[V,e]d^2$. 
So , it contains $4$ or $5$ lines from $\mathcal{L}$. 

Clearly, if $W$ contains $5$ lines from  $\mathcal{L}$,
then these lines from a spread.

If $W$ contains $4$ lines from $\mathcal{L}$, then 
none of these $4$ lines is in $C_V(d)$. So, by \Ref{intersection},
these $4$ lines together with $C_V(d)$ form a spread.
\end{proof}

\bigskip

A spread is called {\em full} if it contains $5$ lines from $\mathcal{L}$.
The set of all full spreads is denoted by $\mathcal{F}$.
A subspace $W$ as in \Ref{defspread} containing exactly $4$ lines from $\mathcal{L}$
contains a unique spread containing these four lines. This spread is called a {\em tangent} spread. By $\mathcal{T}$ we denote the set of tangent spreads.
The fifth $2$-space in such a tangent spread is called the {\em singular} line of the spread.
The set of all singular lines is denoted by $\mathcal{L}_S$.

If we identify each  (full or tangent) spread with the set of lines
from $\mathcal{L}$ contained in it, then  
$(\mathcal{L},\mathcal{F}\cup\mathcal{T})$ is a partial linear space, which we denote by
$\Delta$.


\begin{nummer}\label{2-trans}
Let $W$ be a $4$-space containing a full or tangent spread $S$.
If $S$ is a tangent   spread, then $\langle D_W\rangle$ induces
$\A_4$
on the lines of $\mathcal{L}$ in $S$.
If $S$ is full,  then $\langle D_W\rangle$ induces
$\A_5$
on the lines of $\mathcal{L}$ in $S$.
\end{nummer}

\begin{proof}
Let $l$ be a line of $S$.
An element $d\in D_l$ induces a $3$-cycle on the lines in $S$ and fixes
$l$. 
So, if $S$ is a tangent spread, then  $\langle D_S\rangle$ induces
the $2$-transitive group $\A_4$ on the $4$ lines in $S$.

If $S$ is full, then let $m$ be the unique line of $S$ different from $l$
fixed by $d$.
If there is an element $e\in D_m$ not fixing $l$, then $\langle d,e\rangle$
induces the $2$-transitive group $\A_5$ on $S$.

So, assume that $l\subseteq A_m$, then any line $k\neq l,m$ in $S$
is not in $A_k$, see \Ref{sym}. Let $k$ be such a line and $f\in D_k$ an element
not fixing $m$. Then $\langle d,e,f\rangle$ induces the 
the $2$-transitive group $\A_5$ on $S$.
\end{proof}

\begin{nummer}
If $S$ is a full or tangent spread and $d\in D$, then $C_V(d)$ either contains $S$ or
$C_V(d)\cap S$ is 
a line in $\mathcal{L}\cup\mathcal{L}_S$. 
\end{nummer}
 
\begin{proof}
This follows from \Ref{intersection}.
\end{proof}

\begin{nummer}\label{LcapH}
If $l\in \mathcal{L}$ and $h\in \mathcal{L}_S$, then $h\cap l=\{0\}$.
\end{nummer}

\begin{proof}
Let $S$ be a tangent spread containing $h$ and let $W$ be the $4$-space containing $S$.
Fix an element  $d\in D$ with $l=[V,d]$.
Suppose $l$ meets $h$ nontrivially.
Then $C_V(d)\cap W$ is a line $m\in\mathcal{L}$ of $S$.
Now let $e\in D$ with $[V,e]$ a line of $S$ distinct from $m$.
Then $[V,e]+[V,d]$ contains a spread $T$ and meets $W$ in $[V,e]+(l\cap h)$.
A point in $[V,e]+(l\cap h)$ not on $[V,e]$ or $h$ is on a line 
of $\mathcal{L}$ inside $S$ and on some line
in   $\mathcal{L}\cup \mathcal{L}_S$ of $T$.
Since there is at most one line of  $\mathcal{L}_S$ in $T$,
there is a point in  $[V,e]+(l\cap h)$ on two distinct
lines of $\mathcal{L}$, which contradicts $\mathcal{P}$ being empty. 
\end{proof}

\begin{nummer}\label{inv}
Suppose $W$ is a $4$-space containing  a tangent spread
with singular line $h$.
Suppose  $f\in D$ does not centralize $h$.

Then there exists an involution $t\in \langle D_W\rangle$  with
$h=[V,t]$ being the singular line of the spread, that does not centralize $[V,f]$.
\end{nummer}

\begin{proof}
Let $S$ be the tangent spread in $W$ and $h$ its singular line.
Let $f\in D$ be an element not centralizing $h$. Then, $f$ 
centralizes a line, say $[V,d]$ of $S$, with $d\in D$.
If $d$ does not centralize $[V,f]$, then $[V,d]+[V,f]$
contains a full spread. So, after replacing $d$ by a suitable element in
$D_{[V,d]}$, we can assume that $d$ centralizes $[V,f]$, see \Ref{sym}.
Let $e$ be an element in $D$ with $[V,e]$ a line in $S$ different from 
$[V,d]$. If $e$ centralizes $[V,f]$, then $[V,f]+[V,e]$
is a full spread. So again, after replacing $e$ by a suitable element in
$D_{[V,e]}$, we can assume that $e$ does not centralize $[V,f]$.

The group $\langle d,e\rangle$ induces $\A_4$ on $S$.
After replacing $e$ by $e^{-1}$, if necessary, the element $g=de$
induces a product of two disjoint $2$-cycles in $\A_4$ on $S$.
Then for  
$0\neq v\in [V,f]$, we have 
$0\neq [v,g]=[v,de]=[v,e]\in [V,e]$
and $[v,g^2]=[[v,g],g]\neq 0$ as 
$[V,e]\cap C_V(g)=\{0\}$. 

Let $t$ be the element  $g^2\neq 1$. 
Then $t$ does not centralize $[V,f]$.
The element $t$ fixes all five lines of $S$.
In particular, it fixes $h$ pointwise. 
Any $3$-space on a line $l\in S\setminus \{h\}$
meets all  other lines form $S$, including $h$, in exactly one point.
So, such $3$-space is not only $t$-invariant, all its points
not on $l$ are fixed by $t$. By varying $l$ and the $3$-space, we find all
points $\mathbb{W}$ to be fixed by $t$.
So, $[V,t]\subseteq W\subseteq C_V(t)$, and we can  deduce 
that $t$ has order $2$. 
Indeed, for all $v\in V$ we have $vt^2+vt=[v,t]t=[v,t]=vt+v$,
and thus $vt^2=v$.
Moreover, $C_V(t)$ contains
$C_V(d)\cap C_V(e)$ and $S$  and therefore is of codimension
at most $2$.
Since $t$ commutes with $g$, we find that $[V,t]$ is contained in $h$.
Indeed, for all $v\in V$ we have
$[[v,t],g]=(vt+v)g+(vt+v)=vgt+vg+vt+v=[vg+v,t]=[[v,g],t]=0$. So
$[V,t]\subseteq C_V(g)\cap S=h$.
Hence, either $t$ is a transvection, or $[V,t]=h$.
If $t$ is a transvection, then there is a line $k$ in $\mathcal{L}$
meeting the axis of $t$ in a point.
But then $k$ and $kt$ meet nontrivially, contradicting $\mathcal{P}$
to be empty. So $t$ is indeed the element we are looking for. 
\end{proof}

\bigskip

\begin{nummer}\label{hyperbolic}
Suppose $h\in\mathcal{L}_S$ and $d\in D$ not centralizing $h$.
Then $[V,d]+h$ is a $4$-dimensional space containing two lines 
from $\mathcal{L}$
and three from $\mathcal{L}_S$, pairwise nonintersecting.
\end{nummer}

\begin{proof}
Fix a tangent spread $S$ containing $h$ and an element $d\in D$
not centralizing $h$.
Let $t$ be an involution as in \Ref{inv}.
Then $[V,d]$ and $[V,d]t$ are two  lines from $\mathcal{L}$  
in $[V,d]+h$, and $h,hd$ and $hd^2$ are three lines from $\mathcal{L}_S$ in $S$.
By construction the three lines in $\mathcal{L}_S$ do not intersect.
So, the result follows by \Ref{LcapH}.
\end{proof}

\bigskip

The five lines from $\mathcal{L}\cup\mathcal{L}_S$ in a subspace
$[V,d]+h$, where  $h\in\mathcal{L}_S$ and $d\in D$ not centralizing $h$,
form a spread of $[V,d]+h$ called a {\em hyperbolic} spread.  
The set of all hyperbolic spreads is denoted by $\mathcal{H}$.

\begin{nummer}\label{HcapH}
If $l, h\in \mathcal{L}_S$ are distinct, then $h\cap l=\{0\}$.
\end{nummer}

\begin{proof}
Suppose $h$ and $l$ meet in a point $p$.
Let $S$ be a tangent spread containing $l$ and suppose $f\in D$ is an
element not centralizing $l$. 
By \Ref{inv} there are involutions $t_l$ and $t_h$  in $G$ with
$[V,t_l]=l$ and $[V,t_h]=h$ not centralizing $[V,f]$.
The space 
 $[V,f]+h$ is a hyperbolic spread meeting
$\mathcal{L}$ in  the lines $[V,f]$ and  $[V,f]t_h$, see \Ref{inv},
where $[V,f]t_h=C_V(f)\cap S$.
Similarly   $[V,f]+l$ is  a hyperbolic spread $T$
meeting $\mathcal{L}$ in the two lines
$[V,f]$ and  $[V,f]t_l$, where $[V,f]t_l=C_V(f)\cap T$.
But since $([V,f]+l)\cap ([V,f]+h)$ is $3$-dimensional,
there is a  point in  $C_V(f)\cap S\cap T$, which has  
to be on two lines from $\mathcal{L}$.
This contradicts our assumption that $\mathcal{P}$ is empty.
\end{proof}

\begin{nummer}\label{intersect}
If $S$ and $T$ are two spreads (full, tangent or hyperbolic), 
then $S\cap T$ is empty or a line
in $\mathcal{L}\cup \mathcal{L}_S$.
\end{nummer}

\begin{proof}
This follows immediately from the assumption that  $\mathcal{P}$ is empty, \Ref{HcapH} and \Ref{LcapH}.
\end{proof}

\begin{nummer}\label{connectedness}
$\Delta$ is connected if and only if $G$ is transitive on $\mathcal{L}$.

If  $\Delta$ is connected, then 
the  diameter of its collinearity graph is at most $2$. 
\end{nummer}

\begin{proof}
Suppose $d,e$ are elements from $D$ 
with $[V,d]$ and $[V,e]$ lines not in a spread. Then $d$ centralizes
$[V,e]$ and $e$ centralizes $[V,d]$.
Since $G=\langle D\rangle$, transitivity on the set $\mathcal{L}$
implies the space
$\Delta$ has to be connected.
 
Suppose $\Delta$ is connected.  Then \Ref{2-trans} implies $G$ to be 
transitive on $\mathcal{L}$.
Now suppose $[V,d], [V,f],[V,g], [V,e]$ is a path of length
$3$ in   the 
collinearity graph of $\Delta$.
Let $S$ be the spread in $[V,f]+[V,g]$.
As both $[V,d]$ and $[V,e]$ are in a spread with at least $3$ lines
inside $S$, there is at least one line in $S$ at distance $1$ from both
$[V,d]$ and $[V,e]$.
So the distance between $[V,d]$ and $[V,e]$ is at most $2$.
This implies that   the  diameter of the  collinearity graph
of $\Delta$ is at most $2$.
\end{proof}

\medskip

From now on assume that $\Delta$ is connected and hence also that
$G$ is tranisive on $\mathcal{L}$.

\begin{nummer}\label{fullthennosingular}
If there exists a full spread, then there are no singular lines.
\end{nummer}

\begin{proof} Suppose $S$ is a full spread and $d\in D$ with $[V,d]\in S$.
If $\mathcal{L}_S$ is nonempty, then by transitivity of $G$ on $\mathcal{L}$
there exists a hyperbolic spread $T$ on $[V,d]$ 
containing a second line $[V,e]$ with $e\in D$ from $\mathcal{L}$
and three lines from $\mathcal{L}_S$.
Let $h$ be a singular line in $T$.
Let $t_h$ be an involution in $G$ as in \Ref{inv}
with $[V,t_h]=h$.
The involution $t_h$ centralizes a line $l$ in $S$ distinct from $[V,d]$,
but it maps $[V,d]$ to $[V,e]$.
So, $l+[V,e]$ is also a full spread.

As the group $\langle D_{[V,e]}\rangle$ is transitive on the $4$ lines in the full spread
$l+[V,e]$, see \Ref{2-trans}, but fixes $[V,d]$, 
we find at least $4$ full spreads on $[V,d]$
inside $S+T$. In particular, there are at least $2+16$ lines from
$\mathcal{L}$ in $S+T$.

Now fix an element $f\in D_l$ not centralizing $[V,d]$,
which exists by \Ref{2-trans}.
Then $f$ centralizes at most one of the three singular lines
of $T$.
So, we find at least $1+3\cdot 2=7$ singular lines in $S+T$.
As no two lines from $\mathcal{L}\cup\mathcal{H}$ intersect, there are at
least
$3\cdot(18+7)=75$ projective points in $S+T$, which
contradicts that $S+T$ has dimension $6$, and thus only $63$ points.
Hence $\mathcal{L}_S$ is empty. 
\end{proof}

\begin{nummer}{\bf Theorem.}\label{F_4linear}
Suppose $\Delta$ is connected and contains a full spread.
Then $\Delta$ is isomorphic to
a projective space of order $4$.
In particular, $G$ preserves an  $\F_4$-structure $V_4$
on $V$. 
Moreover, the group $G$ is isomorphic to ${\rm R}(V_4,\Phi)$ for some
subspace $\Phi$ of $V_4^*$ annihilating $V_4$.
\end{nummer}

\begin{proof}
Let $l,m\in \mathcal{L}$ be distinct lines not in a full spread. 
By \Ref{fullthennosingular}, there is not tangent or hyperbolic spread on
$l$ and $m$ and we have
$l\subseteq A_m$ and $m\subseteq A_l$.

By \Ref{connectedness}
there are two full spreads $S$ on $l$  and $T$ on $m$ meeting
at a line $n$. 
Let $k$ be a line in $S$ distinct from $l$ and $n$. Then $m$ and $k$ are in
a full spread $R$. Inside $R$ we can find a line $h$ which spans a full
spread with $n$ distinct from $S$ and $T$. 
But then none of the lines in $R$ is inside $A_l$ and each of them spans
a full spread together with $l$. So, $l$ is on $5$ full spreads each meeting
$T$
in a line. Hence, at least one of these spreads contains $m$. A contradiction.
Thus any two lines from $\mathcal{L}$ are in a full spread.

Now $(\mathcal{L},\mathcal{F})$ is a linear space  of order $4$.
Moreover, it satisfies the Veblen and Young axiom.
Indeed, suppose $S_1,S_2$ are two spreads on a line $l\in \mathcal{L}$,
and $T_1$ and $T_2$ are two spreads meeting both $S_1$ and $S_2$ at lines
distinct from $l$, then as subspaces of $V$, the intersection $T_1\cap T_2$
is $2$-dimensional and thus, by \Ref{intersect}, a line of $\mathcal{L}$.

Thus $V$ carries a $G$-invariant $\F_4$-structure $V_4$ and we can consider
$G$ to be a subgroup of $\GL(V_4)$.
The elements in $D$ induce reflections on $V_4$.
Each $1$-dimensional subspace of $V_4$ is in $\mathcal{L}$ and thus serves
as center of a reflection.
Let $l\in \mathcal{L}$.
By \Ref{2-trans}, no element from $\mathcal{L}$ is in $A_l$.
As in \cite[6.2]{ccs}, we can conclude that
there   is a subspace $\Phi$ of $V_4^*$
annihilating  $V$ with $G\simeq {\rm R}(V,\Phi)$.  
\end{proof}

\bigskip

From now one we can and do assume that $\mathcal{F}$ is empty.

\begin{nummer}\label{dap}
Suppose $S$ is a tangent spread and $h$ its singular line.
If $d\in D$ does not centralize $h$, then the subspace $W$ generated by $S$ and $[V,d]$
contains $12$ lines and $9$ singular lines.
These $21$ lines and spreads in $W$ form a projective plane.
\end{nummer}

\begin{proof}
Let $S$ be a tangent spread and $h$ its singular line.
Let $d\in D$ not centralizing $h$. Then there is a unique
line $m$ in $S$ centralized by $d$. The space $[V,d]+m$ meets 
$\mathcal{L}$ in just $[V,d]$ and $m$, for otherwise it would be a full spread.
So, on $[V,d]$ there are $3$ tangent and one hyperbolic spread inside $W$.
We now easily deduce that
there are $1+9+2=12$ lines of $\mathcal{L}$ in $W$, each on three
tangent spreads.
Together they form a dual affine plane.

Inside $W$ we find also $9$ singular lines. 
As in the proof of the above theorem, we find that these 
lines in $W$ together with the spreads form a projective plane.
\end{proof}

\bigskip
A set of five singular lines in a $2$-dimensional subspace of $V$ is called a {\em singular} spread.
By $\mathcal{S}$ we denote the set of all singular spreads.

\begin{nummer}\label{ap}
Let   $S$ be a tangent spread  with  singular line $h$. 
If $d\in D$
with $[V,d]$ not in $S$ and $C_V(f)\cap S$ equal to $h$,
then $S+[V,f]$ contains $5$ singular lines contained in a singular spread 
and $16$ lines, together forming a projective plane of order $4$.
\end{nummer}

\begin{proof}
Every line $l$ of $S$ determines a unique tangent spread with $[V,d]$. 
So,  $[V,d]$ is on at least $4$ distinct spreads inside $W:=S+[V,d]$.
Thus there are at least $13$ lines  from $\mathcal{L}$
and $5$ singular lines  from $\mathcal{H}$, which are all in the 
$4$-dimensional space $C_W(d)$. Thus the latter $5$ lines form a singular 
spread $T$.
By similar argument we find that all lines form $\mathcal{L}$ inside $W$ 
are on
$4$ tangent spreads. But that implies that there are $16$ lines from
$\mathcal{L}$ in $W$ forming an affine plane.
The rest follows immediately. 
\end{proof}

\begin{nummer}\label{tangent}
If $d\in D$ centralizes $h\in \mathcal{L}_S$, then there is a tangent 
spread $S$ containing $[V,d]$ and $h$.
\end{nummer}

\begin{proof}
Let $S$ be a tangent spread on $h$. 
If $d$ does not centralize the spread (i.e. is collinear in $\Delta$ 
with some point of $S$), then we are done by \Ref{ap}.
Since the graph $\Delta$ is connected, the result follows.
\end{proof}

\begin{nummer}
If $h,l\in \mathcal{L}_S$ are two singular lines,
then there is a hyperbolic or singular spread containing $h$ and $l$.
\end{nummer}

\begin{proof}
Let $S$ be a tangent spread on $h$ and $d\in D$ with $[V,d]\in S$.
Let $t_h$ be an involution in $G$ with $[V,t_h]=h$ as in \Ref{inv}.
If $d$ does not centralize $l$, then inside the subspace
$[V,d]t_h+S$ we find that $h+l$ contains a hyperbolic spread, see \Ref{dap}.

If $d$ centralizes $h$, then $[V,d]+h$ contains a tangent spread
and we can apply \Ref{ap} to find that $h+l$ contains a singular spread.
\end{proof}

\begin{nummer}\label{singularplane}
If $W$ is a $6$-space in $V$ containing two  spreads, then the
lines from $\mathcal{L}\cup \mathcal{H}$ and spreads in $W$ 
form a projective plane of order $4$.
\end{nummer}

\begin{proof}
Let $S$ and $T$ be two  spreads in $W$ intersecting at a  line
$l$.
As each line $m\in S$ different from $l$ 
forms a spread with each line of $T$, there are at least and hence exactly $21$  lines in $W$. Clearly 
the spreads induce a projective
plane of order $4$ on these $21$ lines. 
\end{proof}

\begin{nummer}{\bf Theorem.}\label{unitary}
Suppose $\Delta$ is connected but does not contain  full spreads.
Then the geometry
$(\mathcal{L}\cup\mathcal{L}_S,\mathcal{T}\cup\mathcal{H}\cup\mathcal{S})$ 
is a projective space
of order $4$. The set $\mathcal{L}_S$ is the set of absolute points of
this projective space with respect to some Hermitian polarity.

In particular, $G$ preserves a nondegenerate 
Hermitian $\F_4$-structure $(V_4,h)$.
Moreover, $G$ is isomorphic to ${\rm FU}(V_4,h)$, with $D$ corresponding to the class
of reflections in $G$.
\end{nummer}

\begin{proof}
By the above we find that the geometry $(\mathcal{L}\cup\mathcal{L}_S,
\mathcal{T}\cup\mathcal{H}\cup\mathcal{S})$ is a linear space.
As in the proof  of \Ref{F_4linear} we can prove this space to be a 
projective space of order $4$.
Thus $V$ carries the structure of a vector space $V_4$ over $\F_4$
invariant under $G$.  In particular, we can consider $G$
to be a subgroup of $\GL(V_4)$.

Let $d\in D$ and $l=[V,d]$, 
then $A_l=C_V(d)$ is a hyperplane of $V_4$.

Now consider a singular line $h$, and let $t_h$ be 
an involution as in \Ref{inv} with $h=[V,t_h]$.
From \Ref{dap} and \Ref{ap}
it is readily
seen that  $A_h:=C_V(t_h)$ is a hyperplane of $V_4$ containing 
precisely those spreads on $h$ that are 
tangent or singular. This shows that the map 
$l\in \mathcal{L}\cup \mathcal{L}_S\mapsto A_l$
is a nondegenerate Hermitian polarity on the projective space $\mathbb{P}(V_4)$.
As each  $d\in D$ induces a unitary reflection with center $[V,d]$ on $V_4$,
the theorem readily follows.  
\end{proof}

\bigskip

Now the Theorems \ref{linear},
\ref{cotriangle}, \ref{F_4linear} and \ref{unitary} certainly imply
Theorem \ref{main}. Actually, they provide a proof for a slightly
more general result, in which the assumption on $D$ being a conjugacy
class can be replaced by $\{[V,d]\mid d\in D\}$ being a $G$-orbit.

\vspace{1cm}

\noindent
\begin{minipage}[t]{7cm}
{\small
Hans Cuypers\\
Department of Mathematics\\
Eindhoven University of Technology\\
P.O. BOX 513\\
5600 MB Eindhoven\\
The Netherlands\\
email: f.g.m.t.cuypers@tue.nl}
\end{minipage}
\end{document}